\documentclass[a4paper,11pt]{amsart}

\usepackage{amssymb}
\usepackage{amsmath}

\theoremstyle{plain}
\newtheorem{theorem}{Theorem}[section]
\newtheorem{corollary}[theorem]{Corollary}
\newtheorem{conjecture}[theorem]{Conjecture}
\newtheorem{lemma}[theorem]{Lemma}
\newtheorem{proposition}[theorem]{Proposition}
\theoremstyle{definition}
\newtheorem{definition}[theorem]{Definition}

\newtheorem{remark}[theorem]{Remark}
\newtheorem{review}[theorem]{Review}

\newcommand{\Z}{\mathbb{Z}}
\newcommand{\Q}{\mathbb{Q}}
\newcommand{\R}{\mathbb{R}}
\newcommand{\C}{\mathbb{C}}
\newcommand{\F}{\mathbb{F}}
\newcommand{\G}{\mathbb{G}}

\newcommand{\pd}{\partial}

\newcommand{\Char}{\operatorname{Char}}
\newcommand{\nrm}{\operatorname{N}}

\newcommand{\Res}{\operatorname{Res}}
\newcommand{\Tor}{\operatorname{Tor}}
\newcommand{\Cor}{\operatorname{Cor}}
\newcommand{\Gal}{\operatorname{Gal}}

\newcommand{\Hom}{\operatorname{Hom}}

\newcommand{\Spec}{\operatorname{Spec}}

\newcommand{\Image}{\operatorname{Im}}
\newcommand{\coker}{\operatorname{coker}}

\newcommand{\Ind}{\operatorname{Ind}}

\newcommand{\divi}{\operatorname{div}}

\begin{document}


\title{Milnor $K$-group attached to a torus \\
and Birch-Tate conjecture}

\dedicatory{To Professor Tatsuo Kimura on the
occasion of his sixtieth birthday}


\author{Takao Yamazaki}
\date{\today}
\address{Mathematical Institute, Tohoku University,
  Aoba, Sendai 980-8578, Japan}
\email{ytakao@math.tohoku.ac.jp}

\begin{abstract}
We formulate (and prove under a certain assumption)
a conjecture relating
the order of Somekawa's Milnor $K$-group attached to a torus $T$
and the value of the Artin $L$-function attached to
the cocharacter group of $T$ (regarded as an Artin representation) 
at $s=-1$.
The case $T=\G_m$ reduces to the classical Birch-Tate conjecture.
\end{abstract}


\maketitle

\section{Introduction}

The Birch-Tate conjecture states that,
for a totally real number field $K$,
the following equality should hold:
\[ |\zeta_K(-1)| = \frac{|K_2(O_K)|}{|W_2(K)|}. \]
Here $\zeta_K(s)$ is the Dedekind zeta function of $K$,
$K_2(O_K)$ is the second $K$-group of
the ring $O_K$ of integers in $K$,
and $W_2(K) = H^0(K, \Q/\Z(2)).$
This equality is proved up to a power of $2$
by Wiles \cite{wiles} .

We shall formulate a conjecture with a coefficients in a torus $T$,
which reduces to the Birch-Tate conjecture recalled above when $T=\G_m$.
Let $T$ be a torus over a number field $K$,
and let $X(T) = \Hom(\G_m, T)$ be the cocharacter group.
Then $X \otimes \C$ is a finite dimensional representation 
of the absolute Galois group $G_K$ of $K$
with finite image, i.e. an Artin representation.
Let $L_K(X(T), s)$ be the Artin $L$-function attached to it.
We set 
$W^T(K) = H^0(K, X(T) \otimes \Q/\Z(2))$.
We write $K^T(K)$ for the Milnor $K$-group
$K(K; T, \G_m)$ attached to $T$ and $\G_m$
introduced by Somekawa \cite{somekawa}.
(We will recall the definition of $K^T(K)$
in \S \ref{section:milnorK}.)
In \S \ref{section:btc},
we define a subgroup $K^T(O_K)$ of $K^T(K)$.
When $T=\G_m$, we have identities
\[ L_K(X(T), s) = \zeta_K(s),~~~ 
  W^T(K) \cong W_2(K),~~~
  K^T(O_K) \cong K_2(O_K).
\]
We propose the following generalization of the Birch-Tate conjecture,
which we shall prove for a certain class of tori.

\begin{conjecture}
Let $T$ be a torus over a totally real number field $K$.
Assume $T$ is split by a totally real field.
Then the equality
\[ |L_K(X(T), -1)| = \frac{|K^T(O_K)|}{|W^T(K)|} \]
should hold.
\end{conjecture}

\begin{remark}
The assumption that $T$ is split by a totally real number field
implies that $L_K(X(T), -1)$ 
is a non-zero (rational) number
(cf. proof of Theorem \ref{main2}).
\end{remark}

In \S 2, we introduce a condition 
for a torus (over an arbitrary field)
to `admit a motivic interpretation',
and prove the following.

\begin{proposition}\label{isom}
A torus split by a meta-cyclic extension
admits a motivic interpretation.
(A finite Galois extension $E/F$ of fields
is called meta-cyclic if all Sylow subgroups
of $\Gal(E/F)$ are cyclic.)
\end{proposition}

We also have some examples of tori
which admits a motivic interpretation
without being split by a meta-cyclic extension
(see Remark \ref{genjac}).
Our main result is the following.

\begin{theorem}\label{main}
Let $L/K$ be an extension of totally real fields,
and let $T$ be a torus over $K$ split by $L$.
If $T$ admits a motivic interpretation,
then the equality
\[ |L_K(X(T), -1)| = \frac{|K^T(O_K)|}{|W^T(K)|} \]
holds up to a power of $2$.
\end{theorem}

This result will be proved in \S 4,
where we also prove an analogous result
for a torus over a global field of positive characteristic.
In \S3, we study $K^T(k)$ 
for a torus $T$ over a local field $k$.

\subsection{Conventions}
For a field $F$,
we fix an algebraic closure $\bar{F}$,
and all algebraic extension of $F$
is supposed to be a subfield of $\bar{F}.$
We write $G_F$ for the absolute Galois group of $F$.
For a torus $T$ over a $F$,
we write $X(T) = \Hom(T, \G_m)$ for the cocharacter group of $T$.

Let $A$ be an abelian group.
For a non-zero integer $n$,
we write $A[n]$ and $A/n$ for the kernel and cokernel
of the map $n: A \to A.$
We define $A_{\Tor} = \cup_n A[n]$ 
(resp. $A_{\divi} = \Image(\Hom(\Q, A) \to \Hom(\Z, A)=A)$)
to be the subgroup of torsion elements
(resp. the maximal divisible subgroup) in $A$.
For a prime number $p$,
we define $A[p^{\infty}] = \cup_n A[p^n]$
(resp. $A_{p-\divi} = \Image(\Hom(\Z[\frac{1}{p}], A) \to \Hom(\Z, A)=A)$) 
to be the subgroup of $p$-primary torsion elements
(resp. the maximal $p$-divisible subgroup) in $A$.
We write $A/\divi = A/A_{\divi}$ and $A/p-\divi = A/A_{p-\divi}$.
When a group $G$ acts on $A$,
we write $A^{G}$ and $A_{G}$ for the invariants and coinvariants
of $A$ by $G$.

\section{Milnor $K$-group attached to a torus}\label{section:milnorK}
In this section, $F$ will be an arbitrary field.

\subsection{Definition and basic properties}
Somekawa \cite{somekawa} has introduced the Milnor $K$-group 
$K(F; G_1,\dots, G_r)$
attached to a family of semi-abelian varieties 
$G_1, \ldots, G_r$ over $F$. 
In this paper, 
we only need a special case
where $G_1=T$ is a torus, $G_2=\G_m$ and $r=2$.
To ease the notation, we put
\[ K^T(F) = K(F; T, \G_m). \]
It is defined as a quotient
\begin{equation}\label{def:milnorK}
  K^T(F) = \bigoplus_{E/F} T(E) \otimes E^* / R, 
\end{equation}
where $E$ runs all finite extensions of $F$,
and $R$ is the group
generated by the elements of the following form:
\begin{itemize}
\item(Projection formula)
Let $E_1/E_2/F$ be a tower of finite extensions,
and let $a \in T(E_1), b \in E_2^*$.
Then
\[ N^{E_1}_{E_2}(a) \otimes b - a \otimes R^{E_1}_{E_2}(b) \]
is a generator of $R$.
Here 
$N^{E_1}_{E_2}: T(E_1) \to T(E_2)$ and
$R^{E_1}_{E_2}: E_2^* \hookrightarrow E_1^*$ are
the norm and restriction maps respectively.
\item(Weil reciprocity)
Let $F(C)$ be a function field of one variable over $F$,
and let $S$ be the set of 
all normalized discrete valuation on $F(C)$ over $F$.
For $v \in S$, we write $O_v$ (resp. $F_v$) 
for the valuation ring (resp. the residue field).
Let $a \in T(F(C))$ and $b, c \in F(C)^*$.
Set $S(b) = \{ v \in S ~|~ v(b) \not= 0\}$.
Assume that $a \in T(O_{v_i})$ if $v \not\in S(b)$.
Then 
\[ \sum_{v \in S(b)} a(v) \otimes \pd_v(b, c) 
  + \sum_{v \in S \setminus S(b)} \pd_v(a, c) \otimes b(v)
\]
is a generator of $R$.
Here $\pd_v$ is the local symbol defined in \cite{serre2},
while $a(v) \in T(F_v)$ and $b(v) \in F_v^*$ denote
the the reduction of $a$ and $b$ respectively.
(Recall that $\pd_v(b, c)$ is the usual tame symbol.)
\end{itemize}
The class of $a \otimes b \in T(E) \otimes E^*$ in $K^T(F)$
is written by $\{a, b\}_{E/F}$.
We recall some properties of this group.

\begin{lemma}[\cite{somekawa} Theorem 1.4.]\label{somekawa-milnor}
The correspondence $\{a, b\}_{E/F} \mapsto N^E_F \{a, b\}$
defines a canonical isomorphism
\[ K^{\G_m}(F) \cong K_2(F), \]
where the right hand side is 
the usual second $K$-group,
and $N^E_F$ is the norm map.
\end{lemma}

We often identify $K^{\G_m}(F)$ with $K_2(F)$
by this isomorphism.

\begin{lemma}[\cite{somekawa} Proposition 1.5]
Let $T$ be a torus over $F$,
and let $n$ be a natural number invertible in $F$.
Then we have a homomorphism
\[ h^T_F: K^T(F)/n \to H^2(F, T[n] \otimes \mu_n) \]
called the Galois symbol.
This map satisfies that
$h^T_F( \{a, b\}_{E/F} ) = \Cor^E_F((a) \cup (b))$
for any finite extension $E/F$, $a \in T(E)$ and $b \in E^*$.
Here $(a) \in H^1(E, T[n])$ 
denotes the image of $a$
by the connecting homomorphism associated to
the exact sequence $1 \to T[n] \to T \overset{n}{\to} T \to 1$,
and similarly for $(b) \in H^1(E, \mu_n)$.
\end{lemma}

\begin{remark}
By the Merkurjev-Suslin theorem \cite{ms},
the Galois symbol $h^T_F$ above
is bijective when $T=\G_m$.
It is conjectured by Somekawa \cite{somekawa} that
$h^T_F$ should be always injective.
In Proposition \ref{prop_to_isom} below,
we shall show the injectivity of $h^T_F$
under a certain assumption on $T$.
However, the surjectivity does not hold in general.
(For example, see Proposition \ref{real}).
See \cite{sy} for a related result.
\end{remark}

\begin{lemma}[\cite{sy} Lemma 3]\label{weilres} 
Let $E/F$ be a finite separable extension of fields.
Let $T$ be a torus over $E$,
and let $S = \Res^E_F T$ be the Weil restriction.
Then, we have an isomorphism
\[  K^T(E) \cong K^{S}(F). \]
\end{lemma}

A sequence of algebraic groups $G' \to G \to G''$ over $F$ 
is called {\it Zariski exact} 
if $G'(E) \to G(E) \to G''(E)$ is exact for any extension $E/F$. 

\begin{lemma}[\cite{sy} Lemma 2]\label{rightexact}
Let $R \to S \to T \to 0$ be 
a Zariski exact sequence of tori over $F$.
Then the sequence
\[
K^R(F) \to K^S(F) \to K^T(F) \to 0
\]
is exact as well.
\end{lemma}

\subsection{Motivic interpretation}\label{sect:motivic}
We recall Lichtenbaum's weight two motivic complex.

\begin{review}\label{rev_lich}
Let $\Z(2)$ be the weight two motivic complex,
which is a two-term complex of discrete $G_F$-modules 
(concentrated on degrees one and two)
constructed by Lichtenbaum \cite{lichtenbaum, lichtenbaum2}.
We recall some properties of $\Z(2)$.
\begin{enumerate}
\item
There is a canonical isomorphism
$H^2(F, \Z(2)) \cong K_2(F)$.
\item 
We have $H^3(F, \Z(2))=0$.
\item
If $n \in \Z$ is invertible in $F$,
then we have a triangle
$\Z(2) \overset{n}{\to} \Z(2) \to \mu_n^{\otimes 2} \to \Z(2)[1]$.
\item
If the characteristic $p$ of $F$ positive,
then we have a triangle
$\Z(2) \overset{p^s}{\to} \Z(2) \to \nu_s(2)[-2] \to \Z(2)[1]$,
where $\nu_s(2)$ is 
the second logarithmic Hodge-Witt sheaf of level $s$.
\item
There is a product map $\Z(1) \otimes^{\mathbb{L}} \Z(1) \to \Z(2)$.
(Here $\Z(1) \cong \G_m[-1]$.)
\end{enumerate}
\end{review}

Let $T$ be a torus over $F$,
and let $X = X(T)$ be the cocharacter group of $T$.
For a finite extension $E/F$,
we have a homomorphism
\begin{equation}\label{def:map}
 T(E) \otimes E^* \cong 
   H^1(E, X \otimes \Z(1)) \otimes H^1(E, \Z(1)) 
  \overset{\cup}{\to} H^2(E, X \otimes \Z(2))
  \overset{N^E_F}{\to} H^2(F, X \otimes \Z(2))
\end{equation}
deduced by the product and norm maps.

\begin{definition}\label{conj:isom}
We say $T$ admits a motivic interpretation
if the homomorphism \eqref{def:map}
induces, via eq. \eqref{def:milnorK}, an isomorphism
\[  K^T(F) \overset{\cong}{\to} H^2(F, X \otimes \Z(2)). \]
\end{definition}

We expects any torus admits a motivic interpretation.
In the next subsection,
we prove this under a certain assumption.

\begin{remark}
\begin{enumerate}
\item
It follows from Lemma \ref{weilres} and Shapiro's lemma that,
if a torus $T$ over $F$ admits a motivic interpretation,
then the base change $T \otimes_F E$ of $T$ 
by a finite separable extension $E/F$
admits a motivic interpretation as well.
\item
In order to prove that
the map \eqref{def:map} factors through $K^T(F)$,
one has to show that it kills $R$ in eq. \eqref{def:milnorK}.
There is no difficulty in proving this for the projection formula.
As for the Weil reciprocity,
it seems that a natural way to prove this is to use
the weight three motivic complex $\Z(3)$
(and to show the vanishing of
$H^3(F(C), X \otimes \Z(3)) 
    \to \oplus_v H^2(F(v), X \otimes \Z(2))
    \overset{\text{sum}}{\to} H^2(F, X \otimes \Z(2))$
where $F(C)$ is the function field 
of an irreducible smooth proper curve $C$ over $F$, 
and $v$ runs all closed points of $C$).
If one used Voevodsky's definition of $\Z(r)$,
this would follow from the Gysin sequence \cite{voevodsky}.
However, Voevodsky's theory is developed 
under the assumption of the resolution of singularity.
Because we will also consider 
the global fields of positive characteristic,
we avoid the use of Voevodsky's theory.
See also \cite{moch} for a related result.
\end{enumerate}
\end{remark}

\subsection{Tori split by a meta-cyclic extension}
We recall some facts from \cite{cts}.
A torus $P$ over $F$ is called quasi-trivial if
$P$ is isomorphic to $\oplus_i \Res^{E_i}_F \G_m$,
where $E_i$ runs a family of finite tensions of $F$.
A torus $Q$ over $F$ is called flasque if
$H^1(E, X(Q))=0$ for all finite extension $E/F$.
A torus $I$ over $F$ is called invertible if
there exists a torus $I'$ over $F$ such that
$I \oplus I'$ is quasi-trivial.
We have implications 
`quasi-trivial $\Rightarrow$ invertible $\Rightarrow$ flasque'.

If $T$ is a torus over $F$ split by $E$,
then there exists an exact sequence
\begin{equation}\label{flasque}
   \qquad 0 \to Q \to P \to T \to 0,
\end{equation}
where $P$ (resp. $Q$) is 
a quasi-trivial (resp. flasque) torus over $F$ split by $E$.
We call \eqref{flasque} a flasque resolution of $T$.
A flasque resolution \eqref{flasque} is unique 
up to a direct summand of a quasi-trivial torus in $P$ and $Q$.

\begin{proposition}\label{prop:isom}
Let $T$ be a torus over $F$,
and let \eqref{flasque} be a flasque resolution of $T$.
If $Q$ is invertible, then $T$ admits a motivic interpretation.
\end{proposition}

\begin{proof}
Review \ref{rev_lich} (1)
and Lemma \ref{somekawa-milnor} show that
a split torus admits a motivic interpretation.
By Lemma \ref{weilres} and Shapiro's lemma,
the same holds for a quasi-trivial torus,
hence also for an invertible torus.

Assume a torus $T$ admits a flasque resolution \eqref{flasque}.
If $Q$ is invertible, 
then $H^1(F', Q)=H^3(F', X(Q) \otimes \Z(2))=0$
for any extension $F'/F$
by Hilbert 90 and Review \ref{rev_lich} (2).
This in particular implies that \eqref{flasque} is Zariski exact,
and we have by Lemma \ref{rightexact}
a commutative diagram with exact rows
\[
\begin{matrix}
  &K^Q(F) & \to & K^P(F) & \to & K^T(F) & \to 0
\\
&\downarrow_{\cong} & & \downarrow_{\cong} & & \downarrow &
\\
&H^2(F, X(Q) \otimes \Z(2)) &\to &
H^2(F, X(P) \otimes \Z(2))  &\to &
H^2(F, X(T) \otimes \Z(2))  &\to 0,
\end{matrix}
\]
showing the well-definedness and bijectivity of 
the right vertical map.
\end{proof}

\begin{proof}[Proof of Proposition \ref{isom}]
It follows from a result of Endo-Miyata \cite{endo-miyata}
(see also \cite{cts}) that
a flasque torus split by a meta-cyclic extension
is always invertible.
Now Proposition \ref{isom} is a consequence of Proposition \ref{prop:isom}.
\end{proof}

\begin{remark}\label{genjac}
We give a few examples of a torus $T$
which satisfies the assumption of Proposition \ref{prop:isom}
without being split by a meta-cyclic extension.
\begin{enumerate}
\item
Let $E/F$ be a finite Galois extension
which is not meta-cyclic.
Let $T$ be the kernel of the norm map $\Res_F^E \G_m \to \G_m$.
Then the dual torus $\check{T}$ of $T$ satisfies
the assumption of Proposition \ref{prop:isom},
since it fits into an exact sequence 
$ 0 \to \G_m \to \Res_F^E \G_m \to \check{T} \to 0$.
\item
Let $C$ be an integral proper curve over $F$
whose normalization is isomorphic to the projective line $\mathbb{P}^1$.
Assume that all singular points on $C$
are of coordinate axes type (cf. \cite{tong}).
Then the generalized Jacobian variety $T$ of $C$
is a torus satisfying  the assumption of Proposition \ref{prop:isom}.
Indeed, such $T$ fits into an exact sequence
\[ 0 \to \oplus_s \Res^{F_s}_F \G_m
     \to \oplus_s \oplus_{t \in S(s)} \Res^{F_t}_F \G_m
     \to T \to 0,
\]
where $s$ runs all singular points of $C$,
and $S(s)$ is the inverse image of $s$ by the normalization map.
\end{enumerate}
\end{remark}

\subsection{A few auxiliary results}

\begin{proposition}\label{prop_to_isom}
Let $T$ be a torus over $F$ and let $X=X(T)$.
Assume that $T$ admits a motivic interpretation,
and that $n \in \Z$ is invertible in $F$.
Then we have an isomorphism
\[
H^0(F, T[n] \otimes \mu_n) \cong H^1(F, X \otimes \Z(2))[n]
\]
and exact sequences
\begin{gather*}
0 \to H^1(F, X \otimes \Z(2))/n
\to H^1(F, T[n] \otimes \mu_n) \to K^T(F)[n] \to 0
\\
0 \to K^T(F)/n
\to H^2(F, T[n] \otimes \mu_n) \to H^3(F, X \otimes \Z(2))[n] \to 0.
\end{gather*}
\end{proposition}

\begin{proof}
This follows from the distinguished triangle
$X \otimes \Z(2) \overset{n}{\to} 
 X \otimes \Z(2) \to T[n] \otimes \mu_n
 \to X \otimes \Z(2)[1]$
deduced from Review \ref{rev_lich} (3).
\end{proof}

\begin{corollary}\label{cor_to_isom}
Let $T$ be a torus over $F$ and let $X=X(T)$.
Assume $T$ admits a motivic interpretation.
Let $p$ be a prime different from the characteristic of $F$.
For a natural number $r$, we put 
$\hat{X}_p(r) = \underset{\leftarrow}{\lim} 
   X \otimes \mu_{p^n}^{\otimes r}.$
Then we have isomorphisms
\[ K^T(F)[p^{\infty}]/\divi
 \cong H^2(F, \hat{X}_p(2))_{\Tor} 
 \cong H^1(F, X \otimes \Q_p/\Q_p(2))/\divi.
\]
(Here $H^2(F, \hat{X}_p(2))$ denotes the continuous Galois cohomology.)
\end{corollary}

\begin{proof}
This proof is almost identical to \cite{tate} Theorem 3.5.
Set $M = T[p] \otimes \mu_p$.
We have a commutative diagram with exact rows
\[
\begin{matrix}
0 \to & 
K^T(F)[p] \to & 
K^T(F) \overset{p}{\to} &
K^T(F) \to &
K^T(F)/p \to &
0
\\
& \uparrow_{\text{surj}} 
& \downarrow_{h} & \downarrow_{h} 
& \downarrow_{\text{inj}} & 
\\
& 
H^1(F, M) \to & 
H^2(F, \hat{X}_p(2)) \to & 
H^2(F, \hat{X}_p(2)) \to & 
H^2(F, M), & 
\end{matrix}
\]
where the left and right vertical arrows are
the maps in Proposition \ref{prop_to_isom},
and $h$ is the `continuous symbol' defined 
by the same way as Tate \cite{tate}.
Since $H^2(F, \hat{X}_p(2))$ has no $p$-divisible subgroup
(cf. \cite{tate} Proposition 2.1),
we see $\ker(h) = K^T(F)_{p-\divi}$ and $\coker(h)_{\Tor} = 0.$
This implies that 
$K^T(F)[p^{\infty}] \to H^2(F, \hat{X}_p(2))_{\Tor}$
is a surjection whose kernel is
$K^T(F)[p^{\infty}]_{\divi} = K^T(F)_{\divi}[p^{\infty}]$.
This proves the first identity.
The second isomorphism is given by \cite{tate} Proposition 2.3.
\end{proof}

\begin{lemma}\label{cdtwo}
Let $T$ be a torus over $F$.
Then $H^3(F, X(T) \otimes \Z(2))$ is 
a torsion group of finite exponent.
Moreover we have $H^3(F, X(T) \otimes \Z(2))[p^{\infty}]=0$
if $cd_p(F) \leq 2$ for a prime $p \not= \Char(F)$,
or if $[F : F^p] \leq p$ for $p = \Char(F)$.
\end{lemma}
\begin{proof}
We set $X(2) = X(T) \otimes \Z(2)$.
We take a finite separable extension $E/F$ which splits $T.$
Then we know $H^3(E, X(2))=0$ by Review \ref{rev_lich} (2).
By the norm argument, 
we see that $H^3(F, X(2))$ is annihilated by $n = [L:K]$.
To prove the second assertion,
we write $n= p^k m$ with $(p, m)=1.$
By Review \ref{rev_lich} (3) (resp. (4) ),
$H^3(F, X(2))[p^{\infty}]$ injects to
$H^3(F, T[p^k] \otimes \mu_{p^k})$
(resp. $H^1(F, X(T) \otimes \nu_k(2))$)
when $p \not= \Char(F)$ (resp. $p = \Char(F)$),
which is trivial by assumption.
\end{proof}

\begin{lemma}\label{divisible}
Let $T$ be a torus over a field $F$ of positive characteristic $p$.
Assume $[F: F^p] \leq p$.
Then $K^T(F)$ and $K^T(F)[p^{\infty}]$ 
are $p$-divisible.
\end{lemma}
\begin{proof}
(Cf. \cite{tate2} p. 205.)
It is enough to show the $p$-divisibility of $K^T(F)$.
We take $x \in T(E), y \in E^*$
where $E/F$ is a finite extension.
Because the norm maps 
$(E^{1/p})^* \to E^*$ and $T(E^{1/p}) \to T(E)$
are bijective,
there exist $x' \in T(E^{1/p}), y' \in (E^{1/p})^*$
such that $N(x')=x, N(y')=y$.
Then we have $\{x, y\}_{E/F} = p \{x', y'\}_{E^{1/p}, F}$,
and we are done.
\end{proof}
%

\section{Local field}
When $k$ is a local field, 
we can prove that any torus over $k$
admits a motivic interpretation.
This is an immediate consequence of Theorem \ref{isom} if $k=\R$
(or $k=\C$).

\subsection{Non-archimedean local field}
\begin{lemma}\label{locallemma}
Let $T$ be a torus over a non-archimedean local field $k$.
\begin{enumerate}
\item
$K^T(k)$ is the direct sum of 
a finite group and a uniquely divisible group.
\item
If $p = \Char(k)>0$,
then $K^T(k)[p^{\infty}]=0.$
\item
Let $k_1/k$ be a finite extension.
Then the norm map $N^{k_1}_k: K^T(k_1) \to K^T(k)$ is surjective.
\item
Let $m$ be a natural number invertible in $k$.
Then, the Galois symbol
\[ K^T(k)/m \to H^2(k, T[m] \otimes \mu_m) \]
is bijective.
\item
Let 
$\hat{X}(r) = \underset{\leftarrow}{\lim} X(T) \otimes \mu_n^{\otimes r}$
and
$\Q/\Z(r)' = \underset{\rightarrow}{\lim} \mu_n^{\otimes r},$
where $n$ runs through natural numbers prime to the characteristic of $k$.
Then we have isomorphisms
\[ K^T(k)_{\Tor} \cong K^T(k)/\divi 
\cong H^2(k, \hat{X}(2)) 
\cong H^1(k, X \otimes \Q/\Z(2)')/\divi
\cong \hat{X}(1)_{G_k}.
\]
\end{enumerate}
\end{lemma}
\begin{proof}
We take a finite Galois extension $k'/k$ that splits $T$.
It is proved in \cite{merkurjev} that 
$K_2(k')$ is the direct sum of 
a finite group and a uniquely divisible group.
By the norm argument,
this shows that $K^T(k)$ is the direct sum of 
a uniquely divisible group and a torsion group of finite exponent.
If $\Char(k)=p>0$, then $K^T(k)[p^{\infty}]$ 
is both divisible (by Lemma \ref{divisible}) 
and of finite exponent, hence trivial.
This proves (2).

We prove (3).
When $\Char(k)=0$,
this is proved in \cite{yama} Proposition 3.1.
The same proof works as well when $p=\Char(k)>0$, 
if $[k_1: k]$ is prime to $p$.
The general case can be reduce to this case.
Indeed, the map 
$K^T(k_1)_{\divi} \to K^T(k)_{\divi}$
induced by $N^{k_1}_k$ is surjective by the norm argument.
Thus it suffice to show the surjectivity of
$K^T(k_1)/{\divi} \to K^T(k)/{\divi}$,
which is equivalent to that of
$K^T(k_1)[n] \to K^T(k)[n]$,
where $n$ is the exponent of $K^T(k)_{\Tor}$.
By (2), we know that $n$ is prime to $p$.
Hence we are reduced to the case $[k_1: k]$ is prime to $p$
by the norm argument.

We prove (4) and (5).
We have a commutative diagram
\[
\begin{matrix}
K^T(k')/m & \cong & 
H^2(k', T[m] \otimes \mu_m) & \cong & T[m]_{G_{k'}}
\\
\downarrow_{\nrm} & & 
\downarrow_{\Cor^{k'}_k} & & \downarrow_{\text{proj.}}
\\
K^T(k)/m & \overset{h}{\to}
& H^2(k, T[m] \otimes \mu_m) & \cong & T[m]_{G_{k}}.
\end{matrix}
\]
The upper horizontal map is bijective 
by the Merkurjev-Suslin Theorem \cite{ms}.
The right vertical map is surjective
because it is induced by the identity map on $T[m]$.
This shows the surjectivity of $h$.
This also shows that the kernel of $h \circ N$ is 
$\sum_{\sigma \in \Gal(k'/k)}(1-\sigma)K^T(k')$,
which is killed by $N$ due to the `projection formula' relation.
In view of the surjectivity of $N$ proved in (3),
this shows (4).
Now (5) is an immediate consequence.

Lastly, we prove (1).
If $n$ is the exponent of $K^T(k)_{\Tor}$,
we have
\[ K^T(k)_{\Tor} \cong K^T(k)/n \cong H^2(k, T[n] \otimes \mu_n),
\]
by (4).
Since the right hand side is a finite group,
we see that
$K^T(k)_{\Tor}$ is finite.
This completes the proof.
\end{proof}

\begin{theorem}\label{isom-local}
Let $T$ be a torus over a non-archimedean local field $k$.
Then $T$ admits a motivic interpretation.
\end{theorem}
\begin{proof}
We take a finite Galois extension $k'/k$
which splits $T$.
We set $X(2) = X(T) \otimes \Z(2).$
We are going to show that 
\eqref{def:map} induces 
the homomorphism $\rho$ fitting into the commutative diagram
\[
\begin{matrix}
K^T(k') & \cong & H^2(k', X(2))
\\
\downarrow_{\nrm} & & \downarrow
\\
K^T(k) & \overset{\rho}{\to} & H^2(k, X(2)).
\end{matrix}
\]
The right vertical map is surjective.
Indeed, setting $T' = \ker[\Res^{k'}_k T \to T]$,
we have a distinguished triangle
\[ X(T') \otimes \Z(2) \to \Res^{k'}_k X(2) \to X(2) 
  \to X(T') \otimes \Z(2)[1],
\]
but we have $H^3(k, X(T') \otimes \Z(2))=0$ by Lemma \ref{cdtwo}.
The left vertical map $N$ is also surjective by
Lemma \ref{locallemma} (3).
Lemma \ref{locallemma} (5) shows that
the kernel of $N$ is generated by the elements of the form
$x - \sigma(x)$ with $x \in K^T(k')$ and $\sigma \in \Gal(k'/k)$.
Such an element is killed in $H^2(k, X \otimes \Z(2))$ as well.
This show the existence and surjectivity of $\rho$.

Since $K^T(k)_{\divi}$ is uniquely divisible,
one sees that $\rho|_{K^T(k)_{\divi}}$ is injective
by the norm argument.
On the other hand, $\rho|_{K^T(k)_{\Tor}}$ is also injective as
the composition
\[ K^T(k)_{\Tor} \cong K^T(k)/n \
\to H^2(k, X(2))/n \to H^2(k, T[n] \otimes \mu_n)
\]
(here $n$ is the exponent of $K^T(k)_{\Tor}$)
is bijective by Lemma \ref{locallemma} (4).
Now the theorem follows from Lemma \ref{locallemma} (1).
\end{proof}

\begin{remark}\label{rem-corank}
If $p$ is a prime different from the residue characteristic of $k$,
then $H^1(k, X \otimes \Q_p/\Z_p(2))_{\divi} = 0.$
If further $T$ has good reduction $T_v$ over the residue field $\F$,
then
$H^1(k, X \otimes \Q_p/\Z_p(2)) \cong T_v(\F)[p^{\infty}]$.
\end{remark}

\subsection{Archimedean local field}
Because $K_2(\C)$ is uniquely divisible,
$K^T(\C)$ is uniquely divisible for any torus $T$ over $\C$.
Any torus $T$ over $\R$ admits a motivic interpretation
by Theorem \ref{isom}.
We see that $K^T(\R)$ is the direct sum
of the finite group $K^T(\R)_{\Tor}$
and the uniquely divisible group $K^T(\R)_{\divi}$.
We need to know the structure of $K^T(\R)_{\Tor}$.
Note that any torus over $\R$ is 
isomorphic to a direct sum of copies of tori 
appearing in the following proposition.

\begin{proposition}\label{real}
We have $K^T(\R)_{\Tor} \cong \Z/2\Z$
(resp. $0$) 
if $T = \G_m$
(resp. if $T = \Res^{\C}_{\R} \G_m$ 
or $\ker[\Res^{\C}_{\R} \G_m \to \G_m]$).
Moreover, for any even natural number $n$, 
the exact sequence
\[ 0 \to K^T(\R)/n \to H^2(\R, T[n] \otimes \mu_n) 
\to H^3(\R, X(T) \otimes \Z(2))
\]
is isomorphic to the following sequence:
\begin{align*}
&0 \to \Z/2\Z \to \Z/2\Z \to 0& 
& \text{if}~T = \G_m 
\\
&0 \to 0 \to 0 \to 0& 
& \text{if}~T = \Res^{\C}_{\R} \G_m 
\\
&0 \to 0 \to \Z/2\Z \to \Z/2\Z& 
& \text{if}~T = \ker[\Res^{\C}_{\R} \G_m \to \G_m].
\end{align*}
\end{proposition}

\begin{proof}
The case $T=\G_m$ is well-known.
The other cases can be deduced from Lemma \ref{weilres}
and the exact sequence 
$1 \to \ker[\Res^{\C}_{\R} \G_m \to \G_m]
\to \Res^{\C}_{\R} \G_m \to \G_m \to 1$.
\end{proof}

\section{Global field}\label{section:btc}
Let $K$ be a global field.
For a place $v$ of $K$, we write $K_v$
for the completion of $K$ with respect to $v$.
For a finite place $v$ of $K$,
we write $\F_v$ (resp. $K_v^{nr}$) for the residue field of $v$
(resp. the maximal unramified extension of $K_v$).
When $K$ is a number field,
we write $O_K$ for the ring of integers in $K$,
and set $C=\Spec(O_K)$.
When $K$ is of positive characteristic,
we assume $K$ is the function field of
a smooth projective irreducible curve $C$ 
over a finite field $\F$.

\subsection{Bloch-Moore exact sequence}\label{blochmoore}
We recall some known results.

\begin{theorem}\label{garland}
\begin{enumerate}
\item (Somekawa \cite{somekawa})
Let $T$ be a torus over $K$.
Set $X = X(T)$ and 
$\hat{X}(r) = \underset{\leftarrow}{\lim} X \otimes \mu_n^{\otimes r}$
where $n$ runs through natural numbers invertible in $K$.
Let $m$ be the order of the finite group $\hat{X}(1)_{G_K}$.
Then we have an exact sequence
\[ K^T(K) \to 
  (\bigoplus_{v \not| \infty} \hat{X}(1)_{G_{K_v}})
  \oplus (\bigoplus_{v | \infty} K^T(K_v)/m)
  \to \hat{X}(1)_{G_K} \to 0.
\]
\item (Moore \cite{moore}, Garland \cite{garland})
There exists an exact sequence
\[ 0 \to WK_2(K) \to K_2(K) \to 
  \underset{v: \text{non complex}}{\oplus} \mu(K_v) 
  \to \mu(K) \to 0,
\]
where $WK_2(K)$ is a finite group (so-called wild kernel).
In particular,
$K_2(K)$ is a torsion group without any $p$-divisible subgroup
for any prime $p$.
In the number field case, 
this also implies the finiteness of
\[ K_2(C) = \ker[K_2(K) \to 
  \underset{v \not| \infty}{\oplus} \F_v^*]. 
\]
(In the function field case, we have $WK_2(K)=K_2(C)$.)
\end{enumerate}
\end{theorem} 
We shall prove the finiteness of
the kernel of the first map in (1)
when $T$ admits a motivic interpretation
in Proposition \ref{finiteness} below.

\subsection{Definition of $K^T(O_K)$ and $K^T(C)$}
Let $T$ be a torus over $K$, and let $X=X(T)$.
By Theorem \ref{garland} and the norm argument,
we see that $K^T(K)$ is a torsion group 
without any $p$-divisible subgroup
for any prime $p$.
Hence, by Lemma \ref{divisible}
we have $K^T(K)[p^{\infty}]=0$ if $p=\Char(K)>0$.

\begin{definition}
For each prime $p \not= \Char(K)$,
we define
\[ K^T(C)[p^{\infty}] 
  = \ker[ K^T(K)[p^{\infty}] 
  \to \underset{v {\not|} p\infty}{\prod} 
     H^1(K_v^{nr}, X \otimes \Q_p/\Z_p(2))^{G_{\F_v}}].
\]
Here the $v$-component of the map is the composition
\[ K^T(K)[p^{\infty}] \to K^T(K_v)[p^{\infty}] 
  \cong H^1(K_v, X \otimes \Q_p/\Z_p(2))
  \to H^1(K_v^{nr}, X \otimes \Q_p/\Z_p(2))^{G_{\F_v}},
\]
where the second isomorphism is given by
Theorem \ref{isom-local} and Remark \ref{rem-corank}.
We then define
\[ K^T(C) = \underset{p \not= \Char(K)}{\oplus} K^T(C)[p^{\infty}]. 
\]
In the number field case, we also write $K^T(C)=K^T(O_K)$.
\end{definition}

\begin{remark}\label{kandcoh}
Let $p$ be a prime different from $\Char(K)$.
\begin{enumerate}
\item
When $T$ admits a motivic interpretation,
we have an isomorphism
$K^T(K)[p^{\infty}] \cong H^1(K, X \otimes \Q_p/\Z_p(2))/\divi$
by Corollary \ref{cor_to_isom}.
The corank of $H^1(K, X \otimes \Q_p/\Z_p(2))$
is $r_2 \dim T$, where $r_2$ is the number of complex places on $K$
by \cite{tate} Corollary to Theorem 6.5.
Hence, if further $K$ is totally real or of positive characteristic, 
then we have
\[ K^T(K)[p^{\infty}] \cong H^1(K, X \otimes \Q_p/\Z_p(2)). \]
\item
When $T$ has good reduction at a finite place $v$, 
we have
\[
H^1(K_v, X \otimes \Q_p/\Z_p(2)) \cong T_v(\F_v)[p^{\infty}]
\]
(cf Remark \ref {rem-corank}).
The map
$K^T(K)[p^{\infty}] \to T_v(\F_v)[p^{\infty}]$
can be interpreted 
by the analogous way as the Hilbert symbol
(cf. \cite{somekawa} \S 3).
In particular, we see $K_2(C) = K^T(C)$ if $T=\G_m$.
\item
Summarizing, 
if $K$ is totally real or of positive characteristic, 
and if $T$ admits a motivic interpretation,
then 
$K^T(O_K)[p^{\infty}]$ is isomorphic to the kernel of
\[ H^1(K, X \otimes \Q_p/\Z_p(2))
  \to \underset{v \not\in S, v \not| p}{\oplus} T_v(\F_v)[p^{\infty}] 
      \underset{v \in S, v \not| p}{\oplus} 
        H^1(K_v^{nr}, X \otimes \Q_p/\Z_p(2))^{G_{\F_v}},
\]
where $S$ is a finite set of places of $K$
including all infinite places and
all places where $T$ has bad reduction.
\end{enumerate}
\end{remark}

\subsection{Hasse principle and the finiteness of $K^T(C)$}
In Proposition \ref{finiteness} below,
we prove the finiteness of $K^T(C)$
when $T$ is a torus which admits a motivic interpretation.
In the proof, we need the following result.

\begin{proposition}[Hasse principle]\label{hasse}
Let $T$ be torus over $K$,
and let $X(2)=X(T) \otimes \Z(2)$.
\begin{enumerate}
\item
For all $i \geq 3$, we have an isomorphism
\[ H^i(K, X(2)) \cong \underset{v | \infty}{\oplus} H^i(K_v, X(2)). \]
\item
Suppose that $T$ admits a motivic interpretation.
Let $L/K$ be a finite Separable extension.
For each infinite place $v$ of $K$,
we choose a place $w(v)$ of $L$ above $v$.
Then we have an isomorphism of finite groups
\[ K^T(K)/N^L_K K^T(L) \cong
\underset{v | \infty}{\oplus} K^T(K_v)/N^{L_{w(v)}}_{K_v} K^T(L_{w(v)}).
\]
\end{enumerate}
(When $\Char(K)>0$, both statements mean that
the left hand sides are trivial.)
\end{proposition}

\begin{remark}
It is possible to compute 
the finite group appearing in (2) explicitly
by using Proposition \ref{real}.
When $T=\G_m$, Proposition \ref{hasse} (2)
is proved in \cite{br, ct}.
See also \cite{yama} Proposition 4.1 for a related result.
\end{remark}

\begin{proof}
Firstly, 
we claim that $H^i(K, X(2))$
is a torsion group of finite exponent for all $i \geq 3$.
This is reduced to the case $T=\G_m$ by the norm argument.
By Review \ref{rev_lich} (2), 
we have $H^3(K, \Z(2))=0$.
We also see $H^4(K, \Z(2))=H^3(K, \Q/\Z(2))$ 
is a torsion group of exponent at most $2$.
For $i \geq 5$, the claim follows from 
the spectral sequence 
$E_2^{m, n}=H^m(K, \mathbb{H}^n(\Z(2)))
 \Rightarrow H^{m+n}(K, \Z(2))$, 
together with the fact that
$\mathbb{H}^n(\Z(2)) = 0$ unless $n=1, 2$.

Let $n_i$ be the exponent of $H^i(K, X(2))$ for $i \geq 3$.
We set $n$ to be the prime to $\Char(K)$-part of $n_i n_{i+1}$
and put $M=T[n] \otimes \mu_n$.
The distinguished triangle
$X(2) \overset{n}{\to} X(2) \to M \to X(2)[1]$
induces a commutative diagram with exact rows:
\[
\begin{matrix}
0 \to & H^i(K, X(2)) & \to & H^i(K, M) &  
\to & H^{i+1}(K, X(2)) & \to 0
\\
& \downarrow_{f_{i}} & & \downarrow_{\cong} & & \downarrow_{f_{i+1}}
\\
0 \to & 
\underset{v|\infty}{\oplus} H^i(K_v, X(2)) & \to & 
\underset{v|\infty}{\oplus} H^i(K, M) &  \to & 
\underset{v|\infty}{\oplus} H^{i+1}(K, X(2)) & \to 0
\end{matrix}
\]
Here the middle vertical map is an isomorphism
by the Poitou-Tate theorem
(cf. \cite{serre} \S 6.3 Th\'eor\`eme B).
This shows the injectivity of $f_i$ for all $i \geq 3$
and (using the injectivity of $f_4$ thus obtained)
the surjectivity of $f_i$ for all $i \geq 3$ as well.
When $\Char(K)=p>0$,
a similar argument using Review \ref{rev_lich} (4)
shows that $H^i(K, X(2))[p^{\infty}]=0$ for all $i \geq 3$.
This completes the proof of (1).

We prove (2).
Let $S$ be the kernel of the norm map $\Res^L_K T \to T$.
We set $Y=X(S)$ and $Y(2)=Y \otimes \Z(2)$.
Then we have a distinguished triangle
$Y(2) \to \Res^L_K X(2) \to X(2) \to Y(2)[1]$.
By the assumption that $T$ admits a motivic interpretation,
this induces the exact sequence 
at the upper row in the following commutative diagram
\[
\begin{matrix}
0 \to & K^T(K)/ K^T(L)& \to & H^3(K, Y(2)) & \to & H^3(L, X(2))
\\
& \downarrow & & \downarrow & & \downarrow
\\
0 \to & 
\underset{v|\infty}{\oplus} N^{L_{w(v)}}_{K_v} K^T(L_{w(v)})& \to & 
\underset{v|\infty}{\oplus} H^3(K_v, Y(2)) & \to & 
\underset{v|\infty}{\oplus} \underset{w|v}{\oplus} H^3(L_{w}, X(2))
\end{matrix}
\]
The exact sequence at at the lower row in the diagram
is deduced in a similar way,
by noting the following facts
($v$ is a place of $K$):
(i) The base change of $\Res^L_K T$ to $K_v$
is isomorphic to $\oplus_{w|v} \Res^{L_w}_{K_v}(T \otimes_K K_v)$.
(ii) The image of the norm map
$N^{L_w}_{K_v}: K^T(L_w) \to K^T(K_v)$
is the same for all $w$ over $v$
(because $L/K$ is separable).
(iii) When $v$ is a finite place,
the norm map $N^{L_w}_{K_v}: K^T(L_w) \to K^T(K_v)$
is surjective by Lemma \ref{locallemma} (3).
Now the assertion follows since
the middle and right vertical maps are bijective by (1).
\end{proof}

\begin{proposition}\label{finiteness}
If $T$ is a torus over $K$ which admits a motivic interpretation,
then $K^T(C)$ is a finite group.
\end{proposition}
\begin{proof}
We take a finite Galois extension $L/K$ which splits $T$.
Let $G = \Gal(L/K).$
For a prime $p \not= \Char(K)$,
we have a commutative diagram with exact rows
\[
\begin{matrix}
& K^T(C_L)[p^{\infty}]_G & \to 
& K^T(L)[p^{\infty}]_G & \to 
& \underset{v \not| p\infty}{\oplus} [\underset{w|v}{\oplus}
    H^1(L_w^{nr}, X \otimes \Q_p/\Z_p(2))^{G_{\F_w}} ]_G
& \to 0
\\
& \downarrow & & \downarrow_{f} & & \downarrow
\\
0 \to
& K^T(C)[p^{\infty}] & \to 
& K^T(K)[p^{\infty}] & \to 
& \underset{v \not| p\infty}{\oplus} 
  H^1(K_v^{nr}, X \otimes \Q_p/\Z_p(2))^{G_{\F_v}} ]
& \to 0.
\end{matrix}
\]
If $T$ has good reduction at $v$,
then the $v$-component of 
the right vertical map is an isomorphism since
it is isomorphic to
\[ (\oplus_{w|v} \hat{X}_p(1)_{G_{L_w}})_G
  \cong \hat{X}_p(1)_{G_{K_v}}.
\]
Hence the kernel of the right vertical map is finite.
By Theorem \ref{garland}, $K^T(C_L)$ is a finite group.
By Proposition \ref{hasse} (2),
the cokernel of $f$ is a finite group
which is trivial if $p \not=2$ or $\Char(K)>0$.
This completes the proof.
\end{proof}

\subsection{Isogeny}
We write $\Q/\Z(2)'=\underset{\rightarrow}{\lim} \mu_n^{\otimes 2}$
where $n$ runs natural numbers prime to $\Char(K)$.
For a torus $T$ over $K$, we set
$W^T(K) = H^0(K, X(T) \otimes \Q_p/\Z_p(2))$.

\begin{proposition}\label{isog}
Let $T_1, T_2$ be tori over $K$ admitting a motivic interpretation.
We assume $T_1$ and $T_2$ are isogenous.
If $\Char K>0$, then the equality
\[ \frac{|K^{T_1}(C)|}{|W^{T_1}(K)|} =
   \frac{|K^{T_2}(C)|}{|W^{T_2}(K)|}
\]
holds.
When $K$ is a number field,
the same equality holds up to a power of $2$,
if $T_1$ and $T_2$ are split by a totally real field.
\end{proposition}
\begin{proof}
We only prove the number field case.
(The function field case is easier.)
Let $S$ be a finite set of places of $K$ including
all infinite places
and all places where $T_1$ or $T_2$ have bad reduction.
We fix an odd prime $p$, and set 
$M_i = X(T_i) \otimes \Q_p/\Z_p(2)$ for $i=1, 2$.
We consider a commutative diagram with exact rows for $i=1, 2$
\[
\begin{matrix}
0 \to &K^{T_i}(O_K)[p^{\infty}]& 
  \to &H^1(K, M_i)&
  \to &\oplus_{v \not| p} H^1(K_v^{nr}, M_i)^{G_{\F_v}}&
  \to 0
\\
& \downarrow & 
& \Vert &
& \downarrow &
\\
0 \to &H^1(O_K[\frac{1}{pS}], M_i)& 
  \to &H^1(K, M_i)&
  \to &\oplus_{v \not| pS} T_{i, v}(\F_v)[p^{\infty}]&
  \to 0,
\end{matrix}
\]
where the lower row is the localization sequence
of the etale cohomology,
and the surjectivity of the upper right horizontal map
is due to Theorem \ref{garland} (1),
which implies the surjectivity of the lower right horizontal map as well.
Thus we have an exact sequence
\[ 0 \to K^{T_i}(O_K)[p^{\infty}]
   \to H^1(O_K[\frac{1}{pS}], M_i)
   \to \underset{v \not| p, v | S}{\oplus}
            H^1(K_v^{nr}, M_i)^{G_{\F_v}}
   \to 0.
\]
The localization sequence also implies that
$W^{T_i}(K) = H^0(O_K[\frac{1}{pS}], M_i)$,
$H^2(O_K[\frac{1}{pS}], M_i) = H^2(K, M_i) = 0$.
(The last group is trivial because
$H^2(K, M_i)$ is a torsion group of finite exponent
(by the norm argument) and 
$H^2(K, M)$ is $p$-divisible for any 
$p$-primary torsion divisible group $M$).

Now we are reduced to showing the following equalities:
\begin{align*}
(i)& ~~
\frac{|H^1(O_K[\frac{1}{pS}], M_1)|}{|H^0(O_K[\frac{1}{pS}], M_1)|}
=
\frac{|H^1(O_K[\frac{1}{pS}], M_2)|}{|H^0(O_K[\frac{1}{pS}], M_2)|},
\\
(ii)& ~~
|H^1(K_v^{nr}, M_1)^{G_{\F_v}}| = |H^1(K_v^{nr}, M_2)^{G_{\F_v}}|.
\end{align*}
The isogeny $f: T_1 \to T_2$ implies an exact sequence
\[ 0 \to C \otimes \mu_n \to M_1 \to M_2 \to 0, \]
where $C = \ker(f)$ 
and $n$ is the $p$-power part of the order of $C.$
Then we have an exact sequence
\[ 0 \to \ker(a) \to H^1(K_v^{nr}, M_1) 
   \overset{a}{\to} H^1(K_v^{nr}, M_2) \to 0,
\]
in which $\ker(a)$ is a quotient of 
a finite group $H^1(\F_v, C \otimes \mu_n)$.
Since
$H^1(\F_v, H^1(K_v^{nr}, M_i)) \cong H^2(K_v, M_i) 
= K^T(K_v) \otimes \Q_p/\Z_p = 0$,
we get an exact sequence
\[ 0 \to \ker(a)^{G_{\F_v}} \to H^1(K_v^{nr}, M_1)^{G_{\F_v}} 
   \to  H^1(K_v^{nr}, M_2)^{G_{\F_v}} 
   \to H^1(\F_v, \ker(a)) \to 0.
\]
Since $\ker(a)$ is a finite $G_{\F_v}$-module,
$(ii)$ follows.

In order to prove $(i)$, it suffice to show 
\[ \frac{|H^0(O_K[\frac{1}{pS}], C \otimes \mu_n)| 
           \cdot |H^2(O_K[\frac{1}{pS}], C \otimes \mu_n)|}
   {|H^1(O_K[\frac{1}{pS}], C \otimes \mu_n)|}
   = 1.
\]
By \cite{neu} Theorem 8.6.14, 
the left hand side is equal to
\[
\underset{v | \infty}{\prod} 
  \frac{|H^0(K_v, C \otimes \mu_n)|}{|| n ||_v}.
\]
By the assumption that both $T_1$ and $T_2$ are
split by a totally real field,
$C \otimes \mu_n$
is isomorphic to a direct sum of copies of $\Z/n\Z$
as $G_{K_v}$-modules for any $v | \infty$.
This completes the proof.
\end{proof}

\subsection{Main result}
We now finish the proof of our main result Theorem \ref{main}.
We recall the statement,
including the function field case.

\begin{theorem}\label{main2}
Let $K$ be a global field,
and let $T$ be a torus over $T$.
Assume that $T$ admits a motivic interpretation.
\begin{enumerate}
\item
Suppose that $K$ is a totally real number field,
and that $T$ is split by a totally real field $L$ over $K$.
Then the equality
\[  |L_K(X(T), -1)| = \frac{|K^T(O_K)|}{|W^T(K)|} \]
holds up to a power of $2$.
\item
Suppose $\Char(K)=p>0$. Then the equality
\[ |L_K(X(T), -1)| = \frac{|K^T(C)|}{|W^T(K)|} \]
holds.
(Both sides are rational numbers prime to $p$.)
\end{enumerate}
\end{theorem}

\begin{proof}
Since $L_K(X(T), s)$ is real analytic
(as $X(T)$ being an integral representation),
it is enough to show the equality
after taking the $m$-th power for some $m \in \Z_{>0}$.
Both sides of the equation
is stable under isogeny (by Proposition \ref{isog}).
By \cite{ono}, 
there exist tori $P, Q$ over $K$ such that
\begin{itemize}
\item
Both $P$ and $Q$ are quasi-trivial and split by $L/K$.
\item
$T^{\oplus m} \oplus P$ is isogenous to $Q$
for some $m \in \Z_{>0}$.
\end{itemize}
Hence we are reduced to the case $T=\Res^{M}_K \G_m$ 
for a subextension $M/K$ of $L/K$.
In this case, we have
$L_K(X(T), s) = L_K(\Ind^{M}_K X(\G_m), s) 
= \zeta_M(s)$,
$K^T(C) =K_2(C_{M})$
and $W^T(K) = W_2(M)$.
Thus we are reduced to the case $T=\G_m$, 
which is a theorem of Wiles \cite{wiles} in the number field case,
or of Tate \cite{tate2} in the function field case.
\end{proof}

\vspace{3mm}
\noindent
{\it Acknowledgement.}
Most part of this work was done 
while the author stayed at Universit{\"a}t Bielefeld 
supported by SFB 701. 
He would like to thank Michael Spiess, Bruno Kahn and Joost van Hamel
for stimulating discussion.
This paper is dedicated to Professor Tatsuo Kimura,
whom the author wishes to express his sincerest gratitude
for warm encouragement.



\end{document}